\documentclass[12pt]{article}

\usepackage{amscd, amsmath, amssymb, fancyhdr, epsfig, url}
\usepackage{color}
\numberwithin{equation}{section}

\newcommand{\version}{revised version,\ \   January 13, 2014}

\def\eqref#1{(\ref{#1})}

\newcommand{\arrow}{{\:\longrightarrow\:}}
\newcommand{\Z}{{\Bbb Z}}

\def\1{\sqrt{-1}\:}

\newcommand{\cntrct}                
{\hspace{2pt}\raisebox{1pt}{\text{$\lrcorner$}}\hspace{2pt}}

\makeatletter
\def\x@arrow{\DOTSB\Relbar}
\def\xlongequalsignfill@{\arrowfill@\x@arrow\Relbar\x@arrow}
\newcommand{\xlongequal}[2][]{%
        \ext@arrow 0099\xlongequalsignfill@{#1}{#2}}
\def\xlongrightarrowfill@{\arrowfill@\relbar\relbar\longrightarrow}
\newcommand{\xlongrightarrow}[2][]{%
        \ext@arrow 0099\xlongrightarrowfill@{#1}{#2}}
\makeatother

\renewcommand{\phi}{\varphi}
\renewcommand{\epsilon}{\varepsilon}
\renewcommand{\geq}{\geqslant}
\renewcommand{\leq}{\leqslant}

\newcommand{\Id}{\operatorname{Id}}

\newcounter{Mycounter}[section]
\newcounter{lemma}[section]
\newcounter{claim}[section]
\newcounter{sublemma}[section]
\newcounter{corollary}[section]
\newcounter{theorem}[section]
\renewcommand{\thetheorem}{{Theorem \thesection.\arabic{theorem}}}
\newcommand{\theorem}{%
    \setcounter{theorem}{\value{Mycounter}}
    \refstepcounter{theorem}
    \stepcounter{Mycounter}
    {\noindent \bf \thetheorem.\ }}

\newcounter{conjecture}[section]
\newcounter{proposition}[section]
\newcounter{definition}[section]
\renewcommand{\thedefinition}
      {{Definition~\thesection.\arabic{definition}}}
\newcommand{\definition}{%
    \setcounter{definition}{\value{Mycounter}}
    \refstepcounter{definition}
    \stepcounter{Mycounter}
    {\noindent \bf \thedefinition.\ }}

\newcounter{example}[section]
\newcounter{remark}[section]
\newcounter{problem}[section]
\newcounter{question}[section]
\makeatletter

\renewcommand{\leftmark}{{\scriptsize K\"ahler threefolds without
nontrivial subvarieties}}

\@addtoreset{equation}{section} \@addtoreset{footnote}{section} \makeatother

\def\blacksquare{\hbox{\vrule width 5pt height 5pt depth 0pt}}
\def\endproof{\hfill\blacksquare}

\begin{document}
\begin{center}
{\LARGE\bf
Compact K\"ahler 3-manifolds without\vskip6pt nontrivial subvarieties}\\[4mm]

Fr\'ed\'eric Campana, Jean-Pierre Demailly, Misha Verbitsky\footnote{Misha Verbitsky is
partially supported  by RFBR grants
 12-01-00944-a,  
AG Laboratory NRI-HSE, RF government grant,
ag. 11.G34.31.0023, and NRI-HSE 
Academic Fund Program in 2013-2014, research grant
12-01-0179.\vskip2pt
Submitted on April 29, 2013 to the journal {\em Algebraic Geometry} from  Foundation Compositio Mathematica, accepted in final form on September 4, 2013.
}

\end{center}

{\small \hspace{0.10\linewidth}
\begin{minipage}[t]{0.85\linewidth}
{\bf Abstract} \\
We prove that any compact K\"ahler 3-dimensional manifold
which has no nontrivial complex subvarieties is a torus. This is a very special case of a general conjecture on the structure of so-called simple manifolds, central in the bimeromorphic classification of compact K\"ahler manifolds. The proof follows from the Brunella pseudo-effectivity theorem, combined with fundamental results of Siu and of the second author on the Lelong numbers of closed positive $(1,1)$-currents, and with a version of the hard Lefschetz theorem for pseudo-effective line bundles, due to Takegoshi and Demailly-Peternell-Schneider. In a similar vein, we show that a normal compact and K\"ahler 3-dimensional analytic space with terminal singularities and nef canonical bundle is a cyclic quotient of a simple nonprojective torus if it carries no effective divisor. This is a crucial step towards completing the bimeromorphic classification of compact K\"ahler 3-folds.
\end{minipage}}



\section{Introduction}



We consider here connected compact K\"ahler manifolds $X$ of (complex) dimension $n>1$. An irreducible compact analytic subset $Z$ of $X$ will be said to be a {\em subvariety} of $X$. It is said to be {\em nontrivial} if its (complex) dimension is neither $0$ nor $n$.

The bimeromorphic classification of compact K\" ahler manifolds can be reduced, by means of suitable functorial fibrations, to the following two extreme particular cases: either $X$ is projective or $X$ is {\em simple}, which means that any general point $x$ of $X$ is not contained in any nontrivial subvariety of $X$ (a point is {\em general} if it lies in the intersection of countably many dense Zariski open subsets). If $X$ is a nonprojective torus, this new notion of simpleness coincides with the classical one, defined by the absence of nontrivial subtori. The functorial fibrations needed for this reduction are relative algebraic reductions and relative Albanese maps (see {\cite{_Ca80_}, \cite{_Ca85_}, \cite{_F_}) in the category of connected compact K\"ahler manifolds, where the morphisms are the dominant rational maps with connected fibre.

The present text is concerned with the bimeromorphic classification of such simple $X$. Its difficulty, in contrast to the projective case, is not due to the
abundance and complexity of the examples, but to their (expected) scarceness
and simple structure, which makes essentially all usual invariants of the classification vanish.

For example, if $X$ is simple, its algebraic dimension\footnote{Recall that $a(X)\in \{0,...,n\}$ is the transcendence degree of the field of global meromorphic functions on $X$, $a(X)=n$ if and only if $X$ is projective.} $a(X)$ vanishes, which implies, by \cite{_U_}, Lemma 13.1 and Lemma 13.6, that its Albanese map is surjective and has connected fibres. Because the fibres must then have dimension either $0$ or $n$, we get: either $q=0$, or $X$ is bimeromorphic to its Albanese torus. Thus only the case where $q = 0$ needs to be considered.

The known\footnote{It is shown in \cite{_KV1_} that a stable bundle which does not degenerate to a direct sum of stable bundles on the generic deformation of a Hilbert scheme of $K3$ has a compact moduli space; this may lead to new examples of hyperk\"ahler manifolds. The rigidity of the category of coherent sheaves on the general members of these families is established in \cite{_V02_}, \cite{_V03_}.}
examples of simple compact K\"ahler manifolds are (up to bimeromorphic equivalence) general complex tori, $K3$ surfaces, or the general member of the deformation families of hyperk\"ahler manifolds\footnote{Recall that $X$ is hyperk\"ahler it is K\"ahler and admits a holomorphic symplectic form. It is said to be irreducible if, moreover, $\pi_1(X)$ is finite, and $h^{2,0} (X) = 1$.}. `Simple' surfaces are thus classified: they are either tori ($q > 0$), or $K3$ ($q = 0$). But when $n\geq 3$, the classification is open. The situation is however expected to be similar to the surface case in higher dimensions (see Conjecture 1.1), where the only known examples are constructed out of surfaces which are either $K3$, or tori.

\

The following converse was essentially formulated in \cite{_Ca05_}, as Question~1.4.

\

\begin{conjecture}{\em Let $X$ be a simple compact K\"ahler manifold. Then the following hold.
\begin{itemize}
\item[\bf 1.] either $X$ has a finite \'etale cover bimeromorphic to a complex torus, or $H^0(X,\Omega^2_X)$ is generated by some $\sigma$ which is generically symplectic (that is, $n=2m$ is even, and $\bigwedge^m \sigma\neq 0$). 

We should thus have: $\kappa(X)=0$. If $\dim(X)$ is odd, then $X$ should be bimeromorphic to a complex torus, possibly after some finite \'etale cover. Such a cover will be needed, as shown by the Kummer quotient by the~$-1$ involution.

\item[\bf 2.] When $X$ does not contain any nontrivial subvarieties, $X$ should be either a complex torus, or an irreducible hyperk\"ahler manifold\footnote{Notice that a general deformation of a Hilbert scheme of a $K3$ surface has no subvarieties \cite{_V98_}, while a general deformation of a generalized Kummer variety has some subvarieties, partially classified in \cite{_KV2_}.}.

Thus $K_X$ should be trivial in this case, the two cases being distinguished by $q>0$, or $q=0$.

\item[\bf 3.] If $X$ is generically symplectic, then $\pi_1(X)$ should be finite, of cardinality at most $2^m$, where $2m=n$.
\end{itemize}}
\end{conjecture} 

\
Assertion 3 follows from \cite{_Ca85_}, if $\chi({\cal O}_X)\neq 0$ for $X$ simple and generically symplectic. In what follows, we prove the second assertion of the conjecture for $n=3$.

\

This conjecture can be motivated by the conjectural existence of minimal models in the bimeromorphic category of connected compact K\" ahler manifolds. Indeed, if such a theory exists, and if $X$ is simple, we have $a(X)=0$, which brings $\kappa(X)\leq 0$. The possibility $\kappa(X)=-\infty$ is excluded, since $X$ should then be uniruled. This implies that $\kappa(X)=0$.

In this case, a K\"ahler minimal model is a map $X\dasharrow X'$  such that $X'$ has terminal singularities and $K_{X'}\equiv 0$). A (conjectural) Bogomolov-type decomposition (\cite{_bogomolov_}, \cite{_Beauville_c1=0_}) for manifolds with terminal singularities and $K_{X'}\equiv 0$ could be used to represent a finite cover of $X'$ as a product. However $X$ is simple, hence either $X'$ is covered by a simple torus (with finite locus of ramification), or $X'$ carries a symplectic holomorphic $2$-form, establishing the conjecture above.

The above conjecture can thus be reformulated in terms of conjectural minimal models for compact K\" ahler manifolds, defined as normal compact complex spaces bimeromorphic to compact K\"ahler manifolds, and having $\Bbb Q$-factorial terminal singularities. The simple ones should then have a torsion canonical sheaf, and be either of the form $T/G$, quotients of simple nonprojective tori by a finite group of automorphisms, or carry a multiplicatively unique holomorphic $2$-form, symplectic on the regular locus.

\hfill

\theorem\label{_Main_Theorem_}{\em
Let $X$ be a compact K\"ahler 3-dimensional manifold.
Assume that $X$ has no nontrivial subvariety.
Then $X$ is a torus.}
\endproof

\

The proof will be given as \ref{mth}. Although it is an easy combination of known results, it seems worthwile to present it, because it indicates that the above conjecture might be accessible in dimension $3$ with the actually existing techniques. The proof also gives some information on the higher dimensional case (see \ref{hd}).

\

We also show the following result.

\

\begin{theorem}\label{MT'}{\em Let $X$ be a $3$-dimensional K\"ahler normal connected compact complex space with terminal singularities. Assume that $K_X$ is nef, and that $X$ does not contain any effective divisors\footnote{The nonexistence of effective divisors linearly equivalent to a positive multiple of $K_X$ actually suffices for the conclusion.}. Then $X=T/G$ is isomorphic to the quotient of a nonprojective simple torus by a finite cyclic group of automorphisms.}
\end{theorem}

\

We will use the following simple criterion:

\

\begin{lemma}\label{khi}{\em Let $X$ be a simple compact K\"ahler threefold. If $\chi({\cal O}_X)=0$, then $X$ is bimeromorphic to its Albanese torus.}
\end{lemma} 

\

\noindent
{\em Proof.} By the above remarks, it is sufficient to show that $q>0$. But 
$$0=\chi({\cal O}_X)=1-q+h^{2,0}-h^{3,0}.$$
By Kodaira's theorem, $h^{2,0}>0$ since $X$ is not projective, and $h^{3,0}\leq 1$ since $a(X)=0$. Thus $0=1-q+h^{2,0}-h^{3,0}\geq 1-q+1-1=1-q$, and $q>0$. \endproof

\

The proof of \ref{_Main_Theorem_} is thus reduced to showing that $\chi({\cal O}_X)=0$, which is the objective of the following sections.

\

Finally, let us mention that compact (connected) complex manifolds, not necessarily K\"ahler, without nontrivial subvarieties are also quite important in model theory, giving examples of so-called trivial Zariski geometries (\cite{_Moo_}, \cite{_MP_} \cite{_MMT_}). Our arguments below show that compact $3$-dimensional complex manifolds without nontrivial subvarieties are tori, provided their canonical bundle are pseudo-effective. Brunella's theorem and the hard Lefschetz theorem appear to be the only steps of the proof requiring the K\"ahler hypo\-thesis in depth.


\section{Pseudo-effective and nef line bundles}


We refer to \cite{_Demailly:ecole_} for the basic notions and properties of currents on complex manifolds, and only recall briefly some of the facts needed here.

\

Let $X$ be a compact, $n$-dimensional complex manifold and let $L$ be a holomorphic line bundle on $X$. Assume that $L$ is equipped with a singular hermitian metric $h$ with local weights $e^{-\varphi}$ which are locally integrable. Its curvature current is then $\Theta(L,h):=i\partial \overline{\partial }\varphi$, and its first Chern class $c_1(L)$ is represented by the cohomology class $\frac{1}{2\pi}\;[\Theta(L,h)]$ in the Bott-Chern cohomology
group $H_{\rm BC}^{1,1}(X)$, where 
$$H_{\rm BC}^{p,q}(X)=\{\hbox{$d$-closed $(p,q)$-currents}\}/
\{\hbox{$\partial\overline\partial$-exact $(p,q)$-currents}\}.$$
If $X$ is K\"ahler, $H_{\rm BC}^{1,1}(X)$ is
isomorphic to the Dolbeault cohomology group $H^{1,1}(X)$ and can be viewed as
a subspace of $H^2(X,\Bbb C)$.

\

\definition{\em A $(1,1)$ cohomology class $\{\alpha\}\in H_{\rm BC}^{1,1}(X)$ is said to be {\rm pseudo-effective} if it contains a closed positive $(1,1)$-current~$\Theta$. The class $\{\alpha\}$ is said to be {\rm nef} if for a given smooth positive $(1,1)$-form $\omega>0$ on $X$ and every $\varepsilon>0$, it contains a smooth closed $(1,1)$-form $\alpha_\varepsilon$ such that $\alpha_\varepsilon\ge-\varepsilon\omega$ $[$or~alternatively, in the K\"ahler case, if $\{\alpha\}$ is a limit of K\"ahler classes  $\{\alpha+\varepsilon\omega\}$, where $\omega$ is K\"ahler$]$. A class $\{\alpha\}$ is said to be {\rm big} if it contains a K\"ahler current, namely a closed $(1,1)$-current $\Theta$ such that $\Theta\ge\varepsilon
\omega$ for some smooth positive $(1,1)$-form $\omega$ and some $\varepsilon>0$. A line bundle $L$ on $X$ is said to be nef $($resp.\ pseudo-effective, resp.\ big$)$ if its first Chern 
class $c_1(L)\in H_{\rm BC}^{1,1}(X)$ has the same property.}

\

It follows from the Bochner-Kodaira technique that a big line bundle \hbox{$L\to X$} has maximal Kodaira dimension $\kappa(L)=\dim X$. Therefore, a compact complex manifold $X$ carrying a big line bundle is Moishezon, that is, bimeromorphic to a projective algebraic manifold; in particular, $X$ is not simple. Now, we have the following result.

\

\begin{theorem}\label{siu} (\cite{_Demailly_Paun_}){\em Let $X$ be a compact K\"ahler manifold. Let $\{\alpha\}\in H^{1,1}(X)$ be a nef class. If $\int_X\alpha^n>0$ then $\{\alpha\}$ is big. As a consequence, if $X$ carries a nef line bundle $L$ such that $c_1(L)^n>0$, then $L$ is big
and $X$ cannot be simple.}
\end{theorem}

\

This  special case where $\{\alpha\}=c_1(L)$ (the only one that we need here) can also be obtained as a consequence of holomorphic Morse inequalities \cite{_Demailly:Morse_}. It can be seen as a strengthening of the Grauert-Riemenschneider conjecture proved by Siu \cite{_Siu2_} (although the latter is also valid in the non-K\"ahler case).

\

\definition{\em
The {\rm Lelong number} $\nu_x(\Theta(L,h))=\nu(\varphi,x)$ 
of the $(1,1)$-current $\Theta(L,h)$ at $x\in X$, is defined as $\liminf_{z\to x}\frac{\varphi(z)}{\vert z-x\vert}$, for a local metric on $X$ near~$x$.}

\hfill

\definition{\em
 For a positive real number $c>0$, {\rm the Lelong set} $F_c$
 of a (1,1)-current $\eta$ is the set of points 
 $x\in M$ with $\nu(\eta,x) \geq c$.}

\

By a well-known theorem of Siu (\cite{_Siu:1974_}) any Lelong set of a positive, closed current is a complex analytic subvariety of $M$
(the proof of this difficult result has been considerably simplified, using regularization of currents and the Ohsawa-Takegoshi extension theorem, see \cite{_Demailly:Reg_}.)

\

\begin{theorem}\label{nef}(\cite{_Demailly:Reg_}) {\em Let $X$ be a compact complex manifold. Let $L$ be a pseudo-effective holomorphic line bundle on $X$, with singular hermitian metric $h$ with positive curvature current $\Theta(L,h)$. Assume that the Lelong sets of $\Theta(L,h)$ are all zero-dimensional. Then $L$ is nef, and big as soon as there is at least one nonzero Lelong number.}
\end{theorem}

\

\noindent
{\em Proof.} The first assertion is (\cite{_Demailly:Reg_}, Cor.~6.4). The second follows from (\cite{_Demailly:Reg_}, Cor.~7.6). See also \cite{_Verbitsky:parabolic_}, Thm.~3.12.
\endproof

\

\begin{corollary}\label{zeroLelong}{\em If $X$ is a compact K\"ahler manifold without nontrivial subvarieties, and with pseudo-effective canonical bundle, then $K_X$ is nef, and for any singular hermitian metric $h$ on $K_X, m>0$ with positive curvature current, the Lelong number vanish at any point. In particular, the associated multiplier ideal sheaves on $K_X^{\otimes m}$ are all trivial for any $m>0$ $($that is, $e^{-m\,\varphi}$ is integrable for any $m>0)$.}
\end{corollary}

\section{Hard Lefschetz theorem for the cohomology of pseudo-effective line bundles}

We recall the version of the Hard Lefschetz theorem which is going to be used here.

\

\theorem\label{THL}(\cite{_Takegoshi_}, \cite{_Demailly_Peternell_Schneider:ps-eff_}){\em Let $(X, \omega)$ be a compact K\"ahler manifold, of dimension $n$ with K\"ahler form $\omega$, let $K_X$ be its canonical bundle, and let $L$ be a pseudo-effective holomorphic line bundle on $X$ equipped with a singular
Hermitian metric $h$. Assume that the 
curvature $\Theta$ of $(L,h)$ is a positive current on $X$
and denote by ${\cal I}(h)$ the corresponding multiplier ideal sheaf.
Then the wedge multiplication operator $\eta \arrow \omega^i \wedge \eta$
induces a surjective map
\[
H^0(X,\Omega^{n-i}_X\otimes  L\otimes {\cal I}(h))\stackrel {\omega^i \wedge{}\;
{\scriptstyle\bullet~}}{-\kern-8pt-\kern-8pt-\kern-8pt\arrow}H^i(X,K  \otimes  L\otimes {\cal I}(h)).
\]
Here $\omega$ is considered as an element in $H^1(X,\Omega^1_X)$
and multiplication by $\omega$ maps $H^k(X,\Omega^{n-l}_X\otimes F)$
to $H^{k+1}(X,\Omega^{n-l+1}_X\otimes F)$.}
\endproof

\

This theorem was obtained under various hypotheses during the 1990s, see
\cite{_Enoki:semipositive_} and
\cite{_Mourougane_}. The most general form given here is due to \cite{_Demailly_Peternell_Schneider:ps-eff_}, Thm.~2.1.1. It was proved in \cite{_Takegoshi_} when $L$ is nef.

\hfill

\begin{corollary}\label{a=0}{\em Let $X$ be a compact K\"ahler manifold of dimension $n>1$ without any nontrivial subvarieties. Assume that $K_X$ is pseudo-effective. Then $h^i(X,mK_X)\leq{n \choose i}$, for any $i\geq 0$, and the polynomial $P(m):=\chi(X,mK_X)$ is constant, equal to $\chi({\cal O}_X)$.}
\end{corollary}

\

\noindent
{\em Proof.} Since $n>1$ and $X$ has no nontrivial subvarieties, the algebraic dimension of $X$ vanishes: $a(X)=0$. Therefore $h^0(X,E)\leq \hbox{\rm rank}(E)$, for any holomorphic vector bundle $E$ on $X$. Combining with the Hard Lefschetz theorem (\ref{THL}), this gives the first claim of \ref{a=0}. The second claim is clear because a polynomial function $P(m)$ which remains bounded when $m\to+\infty$ is necessarily constant.
\endproof

\

\begin{corollary}\label{n=3}{\em Let $X$ be a compact K\"ahler manifold of dimension $3$ without any nontrivial subvariety. Assume that $K_X$ is pseudo-effective. The  polynomial $P(m):=\chi(X,mK_X)$ is constant, equal to $\chi({\cal O}_X)=0$.}
\end{corollary}

\

\noindent
{\em Proof.} The intersection number $K_X^3$ vanishes (either by \ref{nef}, or because it is the leading term of $P(m)$, up to the factor $3!$). The Riemann-Roch formula therefore gives $P(m):=\frac{(1-12\,m)}{24}\cdot c_1(X)\cdot c_2(X)$. The boundedness of $P(m)$ then implies that $24\cdot\chi({\cal O}_X)=c_1(X)\cdot c_2(X)=0$.\endproof

\

\begin{remark} Arguments close to some of the ones presented here have already been used in the proof of \cite{_Demailly_Peternell_Schneider:ps-eff_}, theorem 2.7.3: if $X$ is a compact K\"ahler manifold with pseudo-effective canonical bundle admitting a metric with weights $\varphi$ having analytic singularities and positive curvature current, then either $H^0(X, \Omega^i_X\otimes{\cal O}(mK_X))\neq 0$ for infinitely many $m>0$ and some $i\geq 0$, or $\chi(X,{\cal O}_X)= 0$.
\end{remark}


\section{Brunella's pseudo-effectivity criterion.}

The rest of our arguments is based on the following strong (and very difficult)
theorem by Brunella.

\

\begin{theorem}\label{Bru}(\cite{_Brunella:pseudoeff_}){\em Let $X$ be a compact K\"ahler manifold with a $1$-dimen\-sional holomorphic foliation $F$ given by a nonzero morphism of vector bundle $L\to T_X$, where $L$ is a line bundle on $X$, and $T_X$ is its holomorphic tangent bundle. If $L^{-1}$ is not pseudo-effective, then the closures of the leaves of $F$ are rational curves $($and $X$ is thus uniruled$)$.}
 \end{theorem}

\

The following corollary had already been observed in \cite{_HPR_}, Prop.~4.2.

\

\begin{corollary}{\em If $X$ is an $n$-dimensional compact K\"ahler manifold with\break $H^0(X,\Omega^{n-1}_X)\neq 0$, then $K_X$ is pseudo-effective if $X$ is not uniruled.}
\end{corollary}

\

\noindent
{\em Proof.} The vector bundle $\Omega^{n-1}_X$ is canonically isomorphic to $K_X\otimes T_X$. Any nonzero section thus provides a nonzero map $K_X^{-1}\to T_X$ and an associated foliation.\endproof

\

\begin{corollary}\label{s+n=3}{\em If $X$ is a $3$-dimensional nonprojective compact K\"ahler manifold, then either $K_X$ is pseudo-effective or $X$ is uniruled. If $X$ is simple, then $K_X$ is pseudo-effective.}
\end{corollary}

\

\noindent
{\em Proof.}  Kodaira's theorem implies
that any compact K\"ahler manifold with $H^{2,0}(X)=0$
is projective. Thus $H^{2,0}(X)\neq 0$, and the preceding corollary applies (since $2=(n-1)$, here) and gives the first claim. 
If $X$ is uniruled, it is not simple, hence the second assertion. \endproof

\

By combining \ref{s+n=3}, \ref{n=3}, and \ref{khi}, we get our main result:

\

\begin{corollary}\label{mth}{\em If $X$ is a $3$-dimensional compact K\"ahler manifold without nontrivial subvarieties, then $X$ is a complex torus.}
\end{corollary}

\

In higher dimensions, we can replace the assumption from \ref{a=0} that $K_X$ was pseudo-effective with the existence of a holomorphic $(n-1)$-form, giving the following result.

\

\begin{corollary}\label{hd}{\em Let $X$ be a compact K\"ahler manifold of dimension $n>1$ without nontrivial subvarieties. Assume that $H^0(X,\Omega^{n-1}_X)\neq 0$. Then we have $h^i(X,mK_X)\leq{n \choose i}$ for any $i\geq 0$ and the polynomial $P(m):=\chi(X,mK_X)$ is constant, equal to $\chi({\cal O}_X)$ $($The intersection numbers $K_X^j\cdot\hbox{\rm Todd}_{n-j}(X)$ thus all vanish when $j>0$, as expected since $K_X$ should then be trivial.$)$}.\end{corollary}

\

The following result (announced as \ref{MT'}) strengthens \ref{_Main_Theorem_}.

\

\begin{theorem}\label{MT"}{\em Let $X$ be a $3$-dimensional K\"ahler normal connected compact complex space with terminal\/\footnote{we refer to \cite{KM} for this notion, which makes sense without change in the analytic category. In particular, the sheaf $K_X$ is then $\Bbb Q$-Cartier.} singularities. Assume that $K_X$ is nef, and that $X$ does not contain any effective divisors\/\footnote{Divisors linearly equivalent to a positive multiple of $K_X$ actually suffice. This second hypothesis should follow from the first if $X$ is `simple'.}. Then $X=T/G$ is isomorphic to the quotient of a nonprojective simple torus by a finite cyclic group of automorphisms.}
\end{theorem}

\

\noindent
{\em Proof.} By \cite{DP03}, Thm.~7.1, $\kappa(X)=0$. Thus $mK_X$ is a trivial line bundle for some $m>0$. Let $X'\to X$ be the cyclic cover of $X$ defined by a nonzero section of $mK_X$ so that $mK_X$ is Cartier. The singularities of $X'$ are thus terminal as well (\cite{KM}, 5.20). They are also rational, since $K_{X'}$ is Cartier (by \cite{KM}, 5.24). We thus have $\chi({\cal O}_{X'})=\chi({\cal O}_{X''})$ for any smooth model $X''$ of $X'$. But 
$$\chi({\cal O}_{X'})=\frac{-K_{X'}.c_2(X')}{24}=0,$$
since $K_{X'}$ is locally free trivial (see for example \cite{Fl}). The conclusion now follows from \ref{khi} \endproof

\

\begin{remark}\label{prefinalrk}  In order to show the conjecture for `simple' K\"ahler threefolds, it would be sufficient:

1.  to show the following statement (which is a form of the abundance conjecture in this context): If $X$ is a simple (normal, terminal, $\Bbb Q$-factorial, K\"ahler) threefold, and $D$ an effective divisor on $X$, then $D$ is not nef. A possible approach might consist in adapting the proof in the projective case. A more direct approach were however desirable.

2. to show the existence of minimal models such as the ones appearing in \ref{MT'} for any compact K\"ahler threefold $X$. $($This step has just been achieved in \cite{HP13}, where a conjecture equivalent to the one above is formulated in dimension $3$.$)$
\end{remark}

\

\begin{remark}\label{finalrk} In \cite{BGL99}, it is shown that the only group $G$ that can act on a $3$-dimensional torus without nontrivial analytic subvarieties is $\Z^2$ acting by~$\pm\Id$. In fact, from their classification (cf.\
the proof of Thm.~4.1 and Prop.~3.6 in \cite{BGL99}), any other nontrivial action corresponds to a decomposable torus or an abelian variety. Therefore $G$ is trivial or $G=\Z^2$ in~\hbox{\ref{MT"}}.
\end{remark}

\

\noindent
{\bf Acknowledgements:} We are grateful to the Tata Institute
of Fundamental research in Mumbai, India, where this paper was conceived, and especially to I. Biswas and A.J. Parameswaran, organizers  of our stay there.


\hfill

\hfill

\hfill

\small{

\noindent
{\sc Fr\'ed\'eric Campana}
{\sc Institut Elie Cartan\\
Universit\'e Henri Poincar\'e\\
BP 239\\
F-54506. Vandoeuvre-les-Nancy C\'edex\\
et: Institut Universitaire de France\\
\tt frederic.campana@univ-lorraine.fr}\\

\noindent
{\sc Jean-Pierre Demailly}\\
{\sc Acad\'emie des Sciences, et:\\
Universit\'e de Grenoble I\\
Institut Fourier, Laboratoire de Math\'ematiques\\
UMR 5582 du CNRS, BP 74\\
100 rue des Maths, F-38402 Saint-Martin d'H\`eres,\\
\tt jean-pierre.demailly@ujf-grenoble.fr}\\

\noindent
{\sc Misha Verbitsky}\\
{\sc  Laboratory of Algebraic Geometry, SU-HSE,\\
7 Vavilova Str. Moscow, Russia, 117312}\\
{\tt verbit@maths.gla.ac.uk, \ \  verbit@mccme.ru\\}

 }

\end{document}